\documentclass[11pt]{article}

\usepackage{amsmath}
\usepackage{amsthm,verbatim,amssymb}
\usepackage{graphics}
\usepackage{graphicx}
\usepackage{subfigure}
\usepackage{bm}

\usepackage[noend]{algpseudocode}

\usepackage{multirow,multicol}
\usepackage{tikz} 
\theoremstyle{plain}

\theoremstyle{definition}

\usepackage{float}
\usepackage{appendix}
\usepackage{chngcntr}
\usepackage{apptools}
\usepackage{lineno}
\usepackage[ruled,vlined]{algorithm2e}
\AtAppendix{\counterwithin{lemma}{section}}
\AtAppendix{\counterwithin{proposition}{section}}
\AtAppendix{\counterwithin{corollary}{section}}

\usepackage{algpseudocode}

\begin{document}
\title{Adversarial neural network methods for topology optimization of eigenvalue problems\thanks{The work was supported in part by the National Natural Science Foundation of China grant 12071149 and the Science and Technology Commission of Shanghai Municipality (No. 22ZR1421900).}}
\author{Xindi Hu \thanks{Shanghai Key Laboratory of Pure Mathematics and Mathematical Practice \& School of Mathematical Sciences, East China Normal University, Shanghai 200241, China. E-mail address: 763212966@qq.com.}
\and Jiaming Weng \thanks{E-mail address: 1462077856@qq.com.}
\and Shengfeng Zhu \thanks{Corresponding author. E-mail: sfzhu@math.ecnu.edu.cn}}

\maketitle

\begin{abstract}
    This research presents a novel method using an adversarial neural network to solve the eigenvalue topology optimization problems.  The study focuses on optimizing the first eigenvalues of second-order elliptic and fourth-order biharmonic operators subject to geometry constraints. These models are usually solved with topology optimization algorithms based on sensitivity analysis, in which it is expensive to repeatedly solve the nonlinear constrained eigenvalue problem with traditional numerical methods such as finite elements or finite differences. In contrast, our method leverages automatic differentiation within the deep learning framework. Furthermore, the adversarial neural networks enable different neural networks to train independently, which improves the training efficiency and achieve satisfactory optimization results. Numerical results are presented to verify effectiveness of the algorithms for maximizing and minimizing  the first eigenvalues.
\end{abstract}
\section{Introduction}
As an optimization problem subject to constraints of geometry and physics, shape design aims to find a shape which optimizes some objective useful in engineering \cite{MPB2003,JS1992}. Shape design for eigenvalue problems has extensive applications \cite{Henrot} including mechanical vibration \cite{CSG2000, CWC2016, Qian2022, Hu2023, MB}, acoustics \cite{ZS2010}, photonic crystal \cite{CSJ2000}, population dynamics \cite{LYY2006}, and more.

In physical constrained optimization, traditional methods based on sensitivity analysis need to use mesh-based numerical discretization methods such as finite element or finite difference either in the optimize-and-discretize \cite{JS1992} or discretize-and-optimize \cite{MPB2003} framework. After discretization, large-scale optimization problems need to be solved, in which high computational costs are required especially in 3D for sufficiently accurate mesh resolutions. 

The aim of this paper is to seek a shape to maximize or minimize the first, i.e., smallest eigenvalue with deep learning and without relying on the classical theoretical framework of shape sensitivity analysis. It is well-known that with the explosion of available data and computational resources, recent advances in machine learning especially deep learning have produced breakthroughs in different disciplines, including image recognition \cite{AKI2012}, cognitive science \cite{HPB2020}, genomics \cite{BAA2015}, etc. Researchers began to explore applying neural networks in scientific computing and engineering problems. The meshless deep neural network of deep learning has made many contributions in solving partial differential equations, e.g., deep Ritz method \cite{WEA2017,WEB2018}. Neural networks have been applied to solve high-dimensional problems \cite{JSKS2018}. Physics-informed neural networks (PINNs) \cite{RMP2019} integrating many physical constraints and establishing continuous-time models and discrete-time models were proposed to solve forward problems and inverse problems. In \cite{CJSZ}, an improved version of Calderón's method using deep convolution neural networks was developed for electrical impedance tomography. Recently, tensor neural networks were developed \cite{WX} for solving eigenvalue problems.

In addition, to better handle the boundary conditions of physical constraints, a form of automatically satisfying the boundary condition was proposed \cite{JK2018} by multiplying a boundary signed distance function on the output of the neural network. For the variational problems which are difficult to find the boundary signed distance function on irregular domains, \cite{JHH2022} proposed an augmented Lagrangian deep learning method. Subsequently, Hendriks et al. \cite{JCA2020} introduced a method for neural networks to automatically satisfy equation constraints, such as the condition of zero divergence, by simply applying the curl operator to the network's output. This highlights that judicious design of neural networks can significantly enhance both training efficiency and outcomes.

In physical and geometric constraint-governed shape optimization domains, neural networks have also been naturally employed. Inspired by the use of PINNs for solving inverse problems, Lu et al.  \cite{LLP2020} introduced a hybrid Physics-Informed Neural Networks (hPINNs), which represent a network architecture automatically satisfying boundary conditions and circumventing the difficulties associated with the penalty methods on boundaries. Furthermore, the density function under a reasonable assumption was constructed and demonstrated the efficiency of hPINNs for shape optimization in  optical holography and Stokes flows. See also the deep learning toolbox for solving partial differential equations \cite{LLX2020}.


In the field of topology optimization, the use of deep learning methods has received considerable attention in recent years. Deng et al.  \cite{Deng2021} introduced a generative design approach that integrates topology optimization and generative models in an iterative fashion to explore new design schemes. Oh et al.\cite{Oh2019} proposed a deep neural network-based level set method for topology optimization parameters, integrating deep neural networks into the traditional level set method to construct an effective structural topology optimization approach. Lei et al. \cite{Lei2019} employed machine learning to directly establish a mapping between optimal structural/topological design parameters and external loads. Furthermore, Huang et al. \cite{Huang2023} introduced a Problem-Independent Machine Learning (PIML) technique to reduce the computational time associated with finite element analysis. These studies represent the application of deep learning methods to topology optimization problems from various angles, highlighting the vast potential of deep learning in the realm of topology optimization.

Recently, Zang et al. \cite{ZYB2020} proposed a weak adversarial network (WAN) method to solve partial differential equations by transforming the variational form into an operator norm minimization problem. Then the weak solution and test function in the weak form are parameterized into the original network and the adversarial network, respectively. In the training stage, the parameters of the two neural networks are alternately updated. Moreover, the WAN method was used to solve the inverse problems \cite{BGY2020}. For a typical shape and topology optimization problem, the density function as a variable to be optimized is related to the physical state. 

The present paper  proposes a deep learning method based on the adversarial neural network for the eigenvalue optimization problems. Different from the traditional mesh-based discretization method, there is no need in the present meshless method to solve repeatly the eigenvalue problem for sensitivity analysis such as explicit gradient computation. It directly used automatic differentiation in machine learning \cite{AGB2015}. In addition, the adversarial neural networks enable different neural networks to be trained independently, which improves the training efficiency and achieve acceptable optimized results.

The rest of the paper is organized as follows. In Section 2, we introduce topology optimization of second-order elliptic and fourth-order biharmonic eigenvalues. Section 3 introduces the structure of the neural network and the training mode, how to handle constraints, and present the loss function of the neural network. Furthermore, the adversarial neural network algorithm is summarized for eigenvalue optimization problems. In section 4, numerical examples are presented to validate the effectiveness of the algorithm. Conclusions follow in Section 5.

\section{Model problems}

Let $\Omega\subset\mathbb{R}^{2}$ $(d=2,3)$ be an open bounded domain with Lipschitz continuous boundary $\partial\Omega$. 
Let $\rho:\Omega\rightarrow\mathbb{R}$ be a discontinuous function:
\begin{equation}
    \rho= \left\{
    \begin{aligned}
    &\rho_1 &&  \, {\rm in} \,  \Omega \backslash \overline{D}, \\ 
    &\rho_2 && \, {\rm in} \, D,
    \end{aligned}
    \right.
\end{equation}
where $D\subset \Omega$ is an unknown domain and $\rho_1$ and $\rho_2$ ($\rho_1\neq \rho_2$) are given positive constants.
Introduce the characteristic function of $D$ as
\begin{equation}\label{Model1}
    \chi_D= \left\{
    \begin{aligned}
    &1 &&  \, {\rm in} \,  D, \\ 
    &0 && \, {\rm in} \, \Omega \backslash \overline{D}.
    \end{aligned}
    \right.
\end{equation}
Let a constant $c\in (\rho_1|\Omega|, \rho_2|\Omega|)$ be prescribed ($\rho_1<\rho_2$), where $|\Omega|$ denotes the Lebesgue measure of $\Omega$.
Define a set of admissible functions 
\begin{equation}\label{A}
 \mathcal{A}=\left\{\rho \Big| \rho=\rho_{1} \chi_{\Omega \backslash D}+\rho_{2} \chi_{D}, \int_{\Omega} \rho \mathrm{d} x=c\right\}.
\end{equation}
Notice that a mass constraint or equivalently a geometry constraint on $D$ has been imposed in \eqref{A}. Let $L^{2}(\Omega)$ be the space of square-integrable functions defined on $\Omega$.
Denote the Hilbert spaces
$$
\begin{aligned}
H^{1}(\Omega):&=\left\{v \in L^{2}(\Omega)\Big| \frac{\partial v}{\partial x_i} \in L^{2}(\Omega),i=1,\cdots,d\right\}, \\
H_{0}^{1}(\Omega):&=\left\{v \in H^{1}(\Omega) \mid v=0  \text { on } \partial \Omega\right\}.
\end{aligned}
$$
\subsection{Second-order elliptic eigenvalue optimization}

A natural application in vibrating structure engineering, e.g., acoustics \cite{LLY,SOF2001,ZS2010}, is to design the material distribution subject to a mass constraint so that the resonant of a drum frequency reaches the maximum or minimum. See similar applications for the eigenvalue optimization of the composite coefficient problem \cite{CSG2000}.

Let us introduce the model problem.
Given a domain $\Omega\subset\mathbb{R}^d\ (d=2,3)$, let $\lambda_1$ be the first eigenvalue of the problem
\begin{equation}\label{composite_equ}
\left\{
\begin{aligned}
-\Delta u+\alpha  \chi_{D} u &=\lambda u && \text { in } \Omega, \\
u &=0 && \text { on } \partial \Omega,
\end{aligned}
\right.
\end{equation}
where $\alpha>0$. Consider the following eigenvalue optimization problem:
\begin{equation}\label{oijwfe}
\begin{aligned}
\min _{\rho\in \mathcal{A}}\lambda_1\quad \text{or} \quad \max _{\rho\in \mathcal{A}}\lambda_1,
\end{aligned}
\end{equation}
where $c\in(0,|\Omega|)$ is prescribed, $\rho_1=0,\rho_2=1$ and thus $\rho=\chi_D$.

Next, we consider the weak formulation of $\eqref{composite_equ}$: Find $(\lambda,u)\in\mathbb{R}^+\times H^1_0(\Omega)$ such that
$$
\int_{\Omega}(\nabla u \cdot \nabla v+\alpha \chi_{D} uv) \mathrm{d} x=\int_{\Omega} \lambda uv \mathrm{d} x \quad \forall v \in H_{0}^{1}(\Omega).
$$
By Rayleigh quotient, we obtain
\begin{equation}\label{lambda_composite}
\lambda_{1}=\min _{0 \neq v \in H^1_0(\Omega)} \frac{\int_{\Omega}|\nabla v|^{2}+\alpha \chi_{D} v^{2} \mathrm{~d} x}{\int_{\Omega} v^{2} \mathrm{~d} x}.
\end{equation}
Thus, the minimization and  maximization of the first eigenvalue \eqref{oijwfe} can be written in the following forms:
\begin{equation}\label{comp_question}
\min\limits_{\rho\in \mathcal{A}}\ \min _{0 \neq v \in H^1_0(\Omega)} \frac{\int_{\Omega}|\nabla v|^{2}+\alpha \rho v^{2} \mathrm{~d} x}{\int_{\Omega}  v^{2} \mathrm{~d} x}  
\end{equation}
and
\begin{equation}\label{comp_question_max}
\max\limits_{\rho\in \mathcal{A}}\ \min _{0 \neq v \in H^1_0(\Omega)} \frac{\int_{\Omega}|\nabla v|^{2}+\alpha \rho v^{2} \mathrm{~d} x}{\int_{\Omega}  v^{2} \mathrm{~d} x},
\end{equation}
respectively.
\subsection{Fourth-order biharmonic eigenvalue optimization}
Consider fourth-order biharmonic eigenvalue optimization with applications in control of plate frequencies. The difference between the two eigenvalue problems is due to different boundary conditions \cite{CWC2016,MB}. We try to find the optimal density distribution in the open bounded domain $\Omega\subset\mathbb{R}^{2}$. Consider the following two eigenvalue problems subject to inhomogeneous clamped and inhomogeneous simply supported conditions, respectively as
\begin{equation}\label{clamped_equ}
\begin{cases}\Delta^{2} u=\lambda \rho u & \text { in } \Omega, \\ 
u=0  & \text { on } \partial \Omega \\
\frac{\partial u}{\partial n}=0 & \text { on } \partial \Omega 
\end{cases}
\end{equation}
and
\begin{equation}\label{supported_equ}
\begin{cases}\Delta^{2} u =\lambda \rho u & \text { in } \Omega, \\
u=0 & \text { on } \partial \Omega \\
\Delta u-(1-\nu)\kappa\frac{\partial u}{\partial n}=0 & \text { on } \partial \Omega,\end{cases}
\end{equation}
where $n$ denotes the unit outward normal, $\kappa$ is the mean curvature, $\nu$ is the Poisson’s ratio satisfying $\left(-1 \leq \nu \leq 0.5 \right)$. We consider to find the optimal $D$ which makes the first eigenvalue of \eqref{clamped_equ} and \eqref{supported_equ} reaches a minimum or maximum \cite{KK,MB}:
\begin{equation}\label{jojwefo}
\min\limits_{\rho\in \mathcal{A}}\ \lambda_{1}  \quad \text{or}  \quad \max\limits_{\rho\in \mathcal{A}}\ \lambda_{1}.
\end{equation}

For the fourth-order problems \eqref{comp_question} and \eqref{comp_question_max}, let us introduce the following Hilbert spaces:
$$
\begin{aligned}
H^{2}(\Omega):&=\left\{v \in H^{1}(\Omega)\Big| \frac{\partial v}{\partial x_1}, \frac{\partial v}{\partial x_2} \in H^{1}(\Omega)\right\}, \\
H_{0}^{2}(\Omega):&=\left\{v \in H^{2}(\Omega)\Big| v=0, \frac{\partial v}{\partial n}=0 \text { on } \partial \Omega\right\}.
\end{aligned}
$$
By the Rayleigh theorem, we obtain respectively for 
\begin{equation}\label{lambda4_bianfen}
\lambda_{1}=\min _{0 \neq v \in H_0^{2}(\Omega)} \frac{\int_{\Omega}(\Delta v)^{2} \mathrm{~d} x}{\int_{\Omega} \rho v^{2} \mathrm{~d} x}
\end{equation}
and
\begin{equation}\label{lambda42_bianfen}
\lambda_{1} = \min _{0 \neq v \in H_{0}^{1}(\Omega) \cap H^{2}(\Omega) } \frac{\int_{\Omega}(\Delta v)^{2} \mathrm{~d} x-\int_{\partial \Omega}(1-\nu) \kappa\left(\frac{\partial v}{\partial n}\right)^{2} \mathrm{~d} s}{\int_{\Omega} \rho v^{2} \mathrm{~d} x}.
\end{equation}
Then, the problem \eqref{jojwefo} can be written as 
\begin{equation}\label{min1_four}
\begin{aligned}
\left( \max\limits_{\rho\in \mathcal{A}} \right)\min\limits_{\rho\in \mathcal{A}}\ \min _{0 \neq v \in H_0^{2}(\Omega)} \frac{\int_{\Omega}(\Delta v)^{2} \mathrm{~d} x}{\int_{\Omega} \rho v^{2} \mathrm{~d} x}\ 
\end{aligned}
\end{equation}
and
\begin{equation}\label{min2_four}
\begin{aligned}
\left( \max\limits_{\rho\in \mathcal{A}} \right)\ \min\limits_{\rho\in \mathcal{A}}\ \min _{0\neq v \in H_{0}^{1}(\Omega) \cap H^{2}(\Omega)} \frac{\int_{\Omega}(\Delta v)^{2} \mathrm{~d} x-\int_{\partial \Omega}(1-\nu) \kappa\left(\frac{\partial v}{\partial n}\right)^{2} \mathrm{~d} s}{\int_{\Omega} \rho v^{2} \mathrm{~d} x},
\end{aligned}
\end{equation}
associated with \eqref{clamped_equ} and \eqref{supported_equ}, respectively.

\section{Neural networks}
Inspired by biology and neuroscience, artificial neural networks are mathematical models with networks between input arrays and output results. An artificial neural network imitates the neuronal network of the human brain to construct artificial neurons and establish the topological connections between them. Such a network can be regarded as a mathematical mapping. We take the feedforward neural network as an example to introduce the components of the neural network and their mathematical formulations \cite{IGY2016}.

\subsection{Network structure}  

\subsubsection{Layer structure of fully connected feedforward neural networks}
The layer structure is the basic component of the network structure. Each layer is composed of multiple neurons. In a fully connected feedforward neural network, all neurons in a layer are connected to neurons in the previous layer. Now we express how two adjacent layers in a fully connected feed-forward neural network are transmitted. Let the $l$-th and the $(l+1)$-th layers of the neural network contain $d_{l}$ and $d_{l+1} $ neurons, respectively. The output of the $l$-th layer denoted by $z_{l}$ is also the input of the $(l+1)$-th layer. Let $\boldsymbol{W}^{(l+1)}$ be the $(l+1)$-th layer's weight matrix consisting of weight vectors for all neurons of that layer. Denote by $\boldsymbol{b}^{(l+1)}$ and $\phi_{l+1}$ the bias vector and the activation function of the $(l+1)$-th layer, respectively. The mathematical relationship of the output vectors between the two consecutive layers reads:
\begin{equation}
\boldsymbol{z}^{(l+1)}=\phi_{l+1}\left(\boldsymbol{W}^{(l+1)} \boldsymbol{z}^{(l)}+\boldsymbol{b}^{(l+1)}\right),
\end{equation}
where $z_{l}\in \mathbb{R}^{d_{l}}$, $z_{l+1}\in \mathbb{R}^{d_{l+1}}$, $\boldsymbol{W}^{(l+1)}\in\mathbb{R}^{d_{l}\times d_{l+1}}$, and $\boldsymbol{b}^{(l+1)}\in \mathbb{R}^{d_{l+1}}$.

It is obvious that the relationship of the output vectors of the two adjacent layers can be summarized as $$\boldsymbol{z}^{(l+1)}=f^{(l+1)}_{d_l,d_{l+1}}\left( \boldsymbol{z}^{(l)}\right),$$
where $f^{(l+1)}_{d_l,d_{l+1}}:\mathbb{R}^{d_{l}}\rightarrow \mathbb{R}^{d_{l+1}}$ denotes the layer mapping from the $l$-th layer with $d_l$ neurons)to the $l+1$-th layer with $d_{l+1}$ neurons. 

\subsubsection{Block structure}
Compared with simple hierarchical stacking, ResNet is a convolutional neural network \cite{KHX2016} characterized by block structure and jump connection. For a simple layered stacked feed-forward neural network, the training gradient will disappear when the number of layers increases to a certain extent. Fortunately, the jump structure solves this problem well through inserting an identity mapping every several layers such that the gradient is not close to zero. ResNet type networks can be regarded as a block structure if we combine the identity map and several layers. Then the whole network can be regarded as multiple learning block-wise connections. 

\begin{figure}[htbp]
\centering
\includegraphics[scale=0.8]{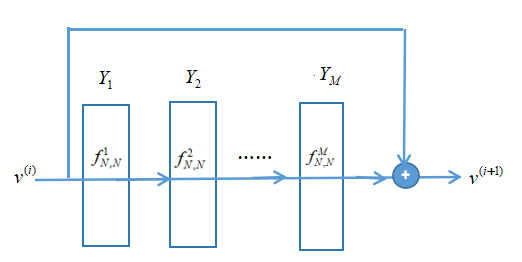}
\caption{Block structure}\label{BS}
\end{figure}
Let us take the jumping hierarchy used in this paper as an example to illustrate the block structure (see Fig. \ref{BS}). This is a block composed of $M$ feed-forward fully-connected layers: $Y_1, Y_2,\cdots, Y_M$ and an identity map. Because the operational dimension of the identity mapping need to be consistent with the output dimension of the $M$-th layer, we simply take the number of neurons in all fully connected layers as equal, which means the number of neuron $d_{l}$ in a layer is degenerated to $N$. If we define $\boldsymbol{v}^{(i)}$ and $\boldsymbol{v}^{(i + 1)}$ as the output of $i$-th and  $(i+1)$-th blocks, then we can express the relationship between $\boldsymbol{v}^{(i)}$ and $\boldsymbol{v}^{(i + 1)}$ :
\begin{equation}
\boldsymbol{v}^{(i+1)}=f^{(M)}_{N,N} \circ...\circ f^{(2)}_{N,N} \circ f^{(1)}_{N,N}\left( \boldsymbol{v}^{(i)}\right)+\boldsymbol{v}^{(i)},
\end{equation}
where $f^{(m)}_{N,N}\, (m=1,2,\cdots,M)$ is the layer mapping consisting of $N$ neurons.
From it, we can find that the relationship of the output of two adjacent blocks can be denoted as $$\boldsymbol{v}^{(i+1)}=g^{(i)}_{M,N}\left( \boldsymbol{v}^{(i)}\right).$$ 
We define $g^{(i)}_{M,N}$ as $i$-th block mapping which consists $M$ layers with $N$ neurons in every layer. 

\subsubsection{Complete structure}
\begin{figure}[ht]
\centering
\includegraphics[scale=0.8]{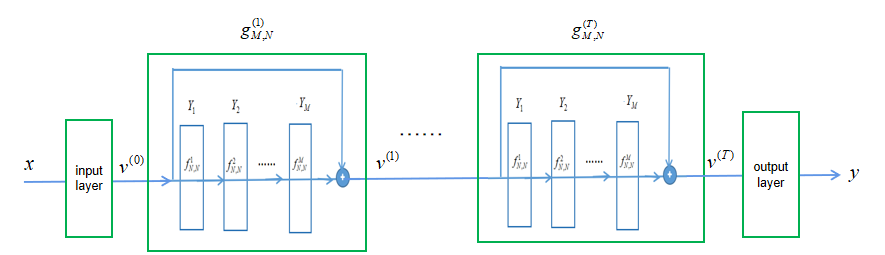}
\caption{Network structure.}\label{structure}
\end{figure}

A complete network structure contains input layers, output layers, and hidden layers or hidden blocks (Fig. \ref{structure}).
According to the above introduction, for an input vector $ \boldsymbol{x} \in \mathbb{R}^{d}$, it will first transform through the input layer, so that it is consistent with the input dimension of the first hidden layer, then into the hidden block, and finally into the output layer, finally getting the output $y$. In our paper, the network is composed of $T$ block hidden blocks. Each block contains the $M$ layers and each layer contains $N$ neurons. Then the transmission process of the whole network is as follows:
\begin{equation}\label{procedure}
\begin{gathered}
\boldsymbol{v}^{(0)}=\boldsymbol{W}_{in}\cdot \boldsymbol{x}, \\
\boldsymbol{v}^{(T)}=g^{(T)}_{M,N}\circ \ldots \circ g^{(1)}_{M,N}(\boldsymbol{v}^{(0)}), \\
y=\phi\big(\boldsymbol{W}_{out} \cdot \boldsymbol{v}^{(T)}+b\big),
\end{gathered}
\end{equation}
where $\boldsymbol{W}_{in}\in \mathbb{R}^{{d}\times{N}}$ and $\boldsymbol{W}_{out}\in \mathbb{R}^{{N}\times{1}}$ are weight matrices of the input layer and output layer, respectively, and $T$ and $N$ are the number of block structures and the number of neurons in the layer structure, respectively. 

Note that the activation function in the output layer can differ from that in the hidden block. Many functions with nonlinear properties can be used as an activation function, which is the key to the neural network dealing with nonlinear problems. Considering the efficiency and stability of neural network training, the activation function generally needs to be continuously differentiable and its derivative can not be too large or too small. In this paper, we mainly use \textbf{tanh function} and use \textbf{square function} when specified.

\subsubsection{Neural network training} 
The previous subsection shows the process of forward transmission of the network. The neural network can be thought as a mapping between an input $ \boldsymbol{x}$ and an output $y$. After building the network, initial parameters require to be set. The Xavier method \cite{GXB2011} is used here to initialize the network. For a specific problem considered to be solved, then a loss function $\mathcal{L}\left(H_{\theta}\right)$ with the neural network $H_{\theta}$ requires to be constructed. Therefore, the training of the neural network is equivalent to the following minimization problem:
$$
\min _{\theta} \mathcal{L}\left(H_{\theta}\right),
$$
where $\theta$ is the set of layer weight matrices and the biases. We try to find the optimal $\theta$ through training the network by backpropagation, which uses the automatic differential method to obtain the derivative \cite{AGB2015}. In the backpropagation optimization method, in addition to the stochastic gradient descent method, there are many optimization algorithms, such as adaptive method, momentum method, Newton method, conjugate gradient method, etc. Here, we use the mainstream method: Adam method \cite{DKJ2016}.

\subsection{Loss functions}\subsubsection{Objective function}

Generally, the loss function needs to consider internal state, boundary and geometry constraints. In this case, the objective function $\eqref{lambda_composite}$ is equivalent to the first eigenvalue pair of equation $\eqref{composite_equ}$. 
Therefor, the loss function only needs to consider the objective $ \eqref{lambda_composite} $, boundary constraint, and geometry constraint.

Firstly, we discrete and approximate the objective based on the equivalence of the objective and $\eqref{composite_equ}$. Instead of using traditional numerical methods such as the finite difference or finite element, we use the PINN method \cite{RMP2019} to build neural networks which we mentioned in Section 2.

We build two adversarial neural networks $H\left(\boldsymbol x  ;\boldsymbol\theta_u\right)$ and $H\left(\boldsymbol x  ;\boldsymbol\theta_\rho\right)$ to approximate $u$ and $\rho$, respectively, where  $\boldsymbol x$ as points in $\Omega$ are inputs of the neural networks, $\boldsymbol\theta_u$ and $\boldsymbol\theta_\rho$ are network parameters that will be updated after initialization. Let $\hat{u}$ and $\hat{\rho}$ be continuously differentiable approximations to ${u}$ and ${\rho}$, respectively. Consider $\hat{u}$ and $\hat{\rho}$ as the outputs of the two neural networks:
$$
\begin{aligned}
\hat{u}\left(\boldsymbol x  ;\boldsymbol\theta_u\right)=H\left(\boldsymbol x  ;\boldsymbol\theta_u\right)\\
\hat{\rho}\left(\boldsymbol x  ;\boldsymbol\theta_{\rho}\right)=H\left(\boldsymbol x  ;\boldsymbol\theta_{\rho}\right).
\end{aligned}
$$

The integrals in the objective $\eqref{lambda_composite}$ can be approximated by numerical quadrature with points in $\Omega$. Let $\left\{x^{(j)},j=1,2,\cdots,N_x\right\}$ be the set of $N_x$ points in $\Omega$.  Then the objective of $\eqref {lambda_composite} $ can be approximated as

\begin{equation*}
\mathcal{L}_{in}\left(\boldsymbol{\theta}_{u}, \boldsymbol{\theta}_{\rho}\right)=
\frac{\sum_{j=1}^{N_x}|\nabla \hat{u}\left(x^{(j)};\boldsymbol{\theta_u}\right)|^{2}
+\alpha\hat{\rho}\left(x^{(j)};\boldsymbol{\theta_\rho}\right)  \hat{u}\left(x^{(j)};\boldsymbol{\theta_u}\right)}{\sum_{j=1}^{N_x}  \hat{u}^{2}\left(x^{(j)};\boldsymbol{\theta_u}\right)}.
\end{equation*}

Now we consider the discretization of the objective for biharmonic eigenvalue optimization.
For the objective $ \eqref{lambda4_bianfen}$ of the clamped plate equation, we approximate it as
\begin{equation}
\mathcal{L}_{in}\left(\boldsymbol{\theta}_{u}, \boldsymbol{\theta}_{\rho}\right)=\frac{\sum_{j=1}^{N_x}|\Delta \hat{u}\left(x^{(j)};\boldsymbol\theta_u\right)|^{2}}{\sum_{j=1}^{N}\hat{\rho}\left(x^{(j)};\boldsymbol\theta_\rho\right)  \hat{u}^{2}\left(x^{(j)};\boldsymbol\theta_u\right)}.
\end{equation}

As for simply supported plate equation, notice that its objective \eqref{lambda42_bianfen} has both volume integral and boundary integral. Introduce $N_b$ points $\left\{x_{b}^{(j)},j=1,2,\cdots,N_{b}\right\}$ on the boundary such that the objective \eqref{lambda42_bianfen} can be approximated as
\begin{equation}\label{obj42_func}
\mathcal{L}_{in}\left(\boldsymbol{\theta}_{u}, \boldsymbol{\theta}_{\rho}\right)=\frac{\frac{A_1}{N_x}\sum_{j=1}^{N_x}|\Delta \hat{u}\left(x^{(j)};\boldsymbol{\theta_u}\right)|^{2}-\frac{A_2}{N_b}\sum_{j=1}^{N_b}(1-\nu)\kappa\frac{\partial \hat{u}}{\partial n}\left(x_b^{(j)};\boldsymbol{\theta_u}\right)^2}{\frac{S_1}{N_x}\sum_{j=1}^{N_x}\hat{\rho}\left(x^{(j)};\boldsymbol{\theta_\rho}\right)  \hat{u}^{2}\left(x^{(j)};\boldsymbol{\theta_u}\right)},
\end{equation}
in which $A_1$ and $A_2$ represent the volume of $\Omega$ and the surface measure of $\partial\Omega$, respectively.

\subsubsection{Boundary constraint}
This section focuses on how to deal with boundary constraints. For neural networks methods in solving partial differential equations (e.g., \cite{WEB2018, RMP2019}), a common method dealing with a boundary constraint is to add a penalty term in objective, which is feasible to many types of boundary conditions: Dirichlet, Neumann, Robin, etc. However, the boundary penalty term probably causes difficulties in training neural networks. In order to satisfy boundary conditions and train the networks easily, we refer to \cite{LLP2020}, in which for the Dirichlet boundary condition, we can introduce a boundary function $\ell_{bd}$ such that $\hat{u}$ directly satisfies the boundary condition automatically. More specifically, set
\begin{equation}\label{hat_u}
\hat{u}\left(\boldsymbol{x} ; \boldsymbol{\theta}_{u}\right)=g(\boldsymbol{x})+\ell_{bd}(\boldsymbol{x}) H\left(\boldsymbol{x} ; \boldsymbol{\theta}_{u}\right)
\end{equation}
to make $\hat{u} $ directly satisfy the Dirichlet boundary condition
\begin{equation}\label{bdy}
u(\boldsymbol{x})=g(\boldsymbol{x}), \quad \boldsymbol{x} \in \partial \Omega.
\end{equation}
In our problems $g(\boldsymbol{x})=0$ and $\ell_{bd}$ satisfies
\begin{equation}\label{jio}
\begin{cases}\ell_{bd}(\boldsymbol{x})=0, & \boldsymbol{x} \in \partial \Omega, \\ \ell_{bd}(\boldsymbol{x})>0, & \boldsymbol{x} \in \Omega .\end{cases}
\end{equation}
Specifically, if the $\partial\Omega$ has a simple geometry structure, an analytical conditional function can be chosen. For example, we can choose $\ell_{bd}=(x-a_1)(x-a_2)(y-a_3)(y-a_4)$ for $\Omega=[a_1,a_2]\times[a_3,a_4]$ with $a_1,a_2, a_3,a_4$ being constants given ($a_1<a_2$ and $a_3<a_4$).

For more complex domains, one can choose a signed distance function \cite{SOR2002} denoted by $f:X\rightarrow\mathbb{R}$ for a larger domain $X\supset\Omega$:
\begin{equation}
f(\bm x)=\left\{
    \begin{aligned}
        &d(\bm x,\partial \Omega) \ &&{\rm if}\,\bm x \in \Omega,\\
       &-d(\bm x, \partial \Omega)\ && {\rm if}\,\bm x \in \Omega^{c},
    \end{aligned}
    \right.
\end{equation}
where the distance function from $\bm x$ to $\partial\Omega$ is defined by
$$
d(\bm x, \partial \Omega):=\inf_{\bm y \in \partial \Omega} d(\bm x, \bm y).
$$
Obviously, the signed distance function  satisfies \eqref{jio} for $ \ell_{bd} $.
For example, when $\Omega= \left\{(x, y): r_{\rm in}^{2}\leq x^{2}+y^{2} \leq\right.$ $\left.r_{\text {out}}^{2}\right\}$, we can choose $\ell_{bd}=\min( r_{\text {out }}^{2}-x^{2}-y^{2},x^{2}+y^{2}-r_{\rm in}^{2} )$, where $r_{\rm in}$ and $r_{\rm out}$ are the inner and outer radiuses, respectively. Notice that the choice of a boundary function $\ell_{bd}$ is not unique.

It is worth mentioning that in \eqref{hat_u}, $\hat{u}$ is no longer the output of the neural network directly, but a continuous differentiable function (with the neural network as a kernel) automatically satisfying the boundary condition \eqref{bdy}.

Next, we deal with the boundary conditions of the fourth-order plate eigenvalue problems. Taking the clamped plate support plate equation as an example, we can find similarly as above a boundary condition function $\ell_{u}$ such that $\hat{u}$ directly meets both the Dirichlet and Neumann boundary conditions. As for $H_0^{2}\left(\Omega\right)$ in this case, we can construct 
\begin{equation}\label{hat_u4}
\hat{u}\left(\boldsymbol{x} ; \boldsymbol{\theta}_{u}\right)=\ell_{u}(\boldsymbol{x}) H\left(\boldsymbol{x} ; \boldsymbol{\theta}_{u}\right),
\end{equation}
where $\ell_{u}$ satisfies
$$
\begin{cases}\ell_{u}(\boldsymbol{x})=0, & \boldsymbol{x} \in \partial\Omega, \\ 
\frac{\partial}{\partial n}\ell_{u}=0, & \boldsymbol{x} \in \partial\Omega,\\
\ell_{u}(\boldsymbol{x})>0, & \boldsymbol{x} \in \Omega,
\end{cases}
$$
with $\frac{\partial}{\partial n}\ell_{u}=\bm n\cdot \nabla \ell_u$.
The disadvantage is that it is sometimes not easy to find the boundary condition function satisfying all above conditions simultaneously. Let us call the above boundary function approach as \textbf {Full-satisfied boundary method}.

We may alternatively use another way to handle the boundary constraints. Taking the simply supported plate problem as an example, since the boundary conditions of the problem are relatively complex, we construct $\hat{u}$ to satisfy the Dirichlet boundary condition, while the other boundary condition considered as a constraint is dealt with  a penalty method
by adding a boundary loss
\begin{equation}\label{bound_func}
\mathcal{L}_b\left(\boldsymbol{\theta}_{u}\right)=M\frac{1}{N_b}
\sum_{j=1}^{N_b} \left[\Delta\hat{u}^2\left(x_b^{(j)};\boldsymbol{\theta_u}\right)^{2}-(1-\nu)\kappa\frac{\partial \hat{u}}{\partial n}\left(x_b^{(j)};\boldsymbol{\theta_u}\right)\right]^2,
\end{equation}
to the total loss function $\mathcal{L}_{sum}$, where a large penalty parameter $M>0$. We call it \textbf{Half-satisfied boundary method}.

\subsubsection{Geometry constraint}
Finally, we deal with the geometry constraint with the penalty method. The original constrained optimization problem is transformed into an unconstrained optimization problem. We first consider numerical quadrature to approximate the geometry constraint. If there are enough, e.g., $N_x$, sampling points in domain $\Omega$, $\int_{\Omega}\rho$ can be approximated as quadrature. Then, the geometry constraint $\int_{\Omega}\rho= c$ can be replaced approximately as $h\left(\boldsymbol{\theta}_{\rho}\right)=0$, where 
$$
h\left(\boldsymbol{\theta}_{\rho}\right)=\frac{\sum_{j=1}^{N_x} \hat{\rho}\left(x_{j};\boldsymbol{\theta}_{\rho}\right)}{N_x}-c.
$$
If the penalty method is used, then the loss function related to the geometry constraint is as follows:
$$
\mathcal{L}_{w}\left(\boldsymbol{\theta}_{\rho}\right)=\mu h^{2}\left(\boldsymbol{\theta}_{\rho}\right),
$$
where a penalty coefficient $\mu>0$ is introduced. In the general penalty method, $\mu$ is a fixed value. To satisfy the constraint better, the fixed value of $\mu$ usually should not be very large. However, a large $\mu$ leads to inefficient training due to suffering very small training rate and thus makes it difficult to reach the global minimizer \cite{DPB2014}. One improvement is to allow $\mu$ to start from a small initial value and increase during iterations \cite{DPB2014}. The loss function of this variable coefficient penalty method is
\begin{equation}\label{mu_k}
\mathcal{L}_{w}\left(\boldsymbol{\theta}_{\rho}\right)=\mu_{k} h^{2}\left(\boldsymbol{\theta}_{\rho}\right),
\end{equation}
where $\mu_{k}>0$ is a penalty coefficient of the $k_\text{th}$ iteration and $\mu_{k+1}=\beta \mu_{k}
$
with $\beta>1$.


Another way compared with the variable coefficient penalty method is the augmented Lagrangian method \cite{TFC2003} with
\begin{equation}\label{lambda_k}
\mathcal{L}_{w}\left(\boldsymbol{\theta}_{\rho}\right)=\mu_{k}h^{2}\left(\boldsymbol{\theta}_{\rho}\right)+\lambda_{k} h\left(\boldsymbol{\theta}_{\rho}\right),
\end{equation}
where the $k$-th Lagrange multiplier $\lambda_{k}$ is updated during iterations:
$$
\lambda_{k}=\lambda_{k-1}+2 \mu_{k-1} h\left( \boldsymbol{\theta}_{\rho}^{k-1}\right).
$$
In addition, $\mu_ {k}$ grows as
$$
\mu_{k+1}=\min(\beta\mu_{k},S),
$$
where $S$ is a upper bound for $\mu_{k}$. A benefit of augmented Lagrangian is the dynamic adjustment, which allows $\mu_{k}$ to reach convergence when $\mu_{k}$ does not need to go to infinity, avoiding slow training caused by too large value of $\mu_{k}$.  
In this paper, we try to treat the geometry constraint with both the variable coefficient penalty method \eqref{mu_k} and the augmented Lagrangian method \eqref{lambda_k}.

\subsubsection{Total loss function}
Now we can write the total loss function $\mathcal{L}_{sum}$. For a minimization problem, we can use \textbf {full-satisfied boundary method} for exactly boundary condition satisfaction, here the total loss function is
\begin{equation}\label{sum4_func_1}
\mathcal{L}_{sum}\left(\boldsymbol{\theta}_{u}, \boldsymbol{\theta}_{\rho}\right)=\mathcal{L}_{in}\left(\boldsymbol{\theta}_{u}, \boldsymbol{\theta}_{\rho}\right)+\mathcal{L}_{w}\left(\boldsymbol{\theta}_{\rho}\right).
\end{equation}
If the boundary conditions cannot be directly satisfied, we use \textbf{half-satisfied boundary method}. Then the total loss $\mathcal{L}_{sum}$ consists of three parts:
\begin{equation}\label{sum4_func_2}
\mathcal{L}_{sum}\left(\boldsymbol{\theta}_{u}, \boldsymbol{\theta}_{\rho}\right)=\mathcal{L}_{in}\left(\boldsymbol{\theta}_{u}, \boldsymbol{\theta}_{\rho}\right)+\mathcal{L}_{w}\left(\boldsymbol{\theta}_{\rho}\right)+\mathcal{L}_b\left(\boldsymbol{\theta}_{u}\right).
\end{equation}

As for maximization problem, since different loss functions can be defined in the inner loop, the adversarial neural network method successfully represent the objective of $\eqref{comp_question_max}$ in the neural network framework. This is an advantage that the general simultaneous neural network method does not have. For this objective of maximum-minimum formulation, when $\hat{u}$ is fixed in inner loop, to find $\hat{\rho}$ that maximizes the first eigenvalue and satisfies the geometry constraint, we only need add a minus sign before $\mathcal{L}_{in}$  in $\mathcal{L}_{sum}$, i.e.
\begin{equation}
\label{sum4_func_3}
\mathcal{L}_{sum}\left(\boldsymbol{\theta}_{u}, \boldsymbol{\theta}_{\rho}\right)=-\mathcal{L}_{in}\left(\boldsymbol{\theta}_{u}, \boldsymbol{\theta}_{\rho}\right)+\mathcal{L}_{w}\left(\boldsymbol{\theta}_{\rho}\right)
\end{equation}
for \textbf {full-satisfied boundary method} and 
\begin{equation}
\label{sum4_func_4}
\mathcal{L}_{sum}\left(\boldsymbol{\theta}_{u}, \boldsymbol{\theta}_{\rho}\right)=-\mathcal{L}_{in}\left(\boldsymbol{\theta}_{u}, \boldsymbol{\theta}_{\rho}\right)+\mathcal{L}_{w}\left(\boldsymbol{\theta}_{\rho}\right)+\mathcal{L}_b\left(\boldsymbol{\theta}_{u}\right)
\end{equation}
for \textbf{half-satisfied boundary method}.

\subsection{Optimization training}
In last section, we have given the unconstrained loss functions under different methods. Now we summarize the complete training process in this section.
 
The first step is to establish two neural networks $H\left(\boldsymbol x  ;\boldsymbol\theta_u\right)$ and $H\left(\boldsymbol x  ;\boldsymbol\theta_\rho\right)$, which are the block structures introduced in section 2. 
Next, we need to select appropriate auxiliary functions to construct $\hat{u}\left(\boldsymbol{x};\boldsymbol{\theta}_{u}\right)$ and $\hat{\rho}\left(\mathbf{x}; \boldsymbol{\theta}_{\rho}\right)$, which both contain neural networks as their main body. To treat boundary constraints mentioned above, we choose an appropriate boundary condition function $\ell_{u}$ attached to $H\left(\boldsymbol x  ;\boldsymbol\theta_u\right)$. Besides, we can also make some appropriate deformation based on $H\left(\boldsymbol x;\boldsymbol\theta_\rho\right)$ so that $\hat{\rho}$ can satisfy some initial guess.
After that, we uniformly sample points in $\Omega$ in preparation for the discretization. 
With the sample points, we can calculate the ${L}_{sum}$ and start to train the neural networks. A common procedure is to update $\boldsymbol{\theta}_{u}$ and $\boldsymbol{\theta}_{\rho}$ simultaneously by minimizing the loss function ${L}_{sum}$ in one iteration. However, this method makes the network parameters of two different physical variables update with the same step size, which reduces the training efficiency. In this paper, we use the adversarial training to improve efficiency. More specifically, the alternating direction method is used to update $\boldsymbol{\theta}_{u}$ and $\boldsymbol{\theta}_{\rho}$.

In one iteration of the inner loop, we first fix $\boldsymbol{\theta}_{\rho}$, update $\boldsymbol{\theta}_{u}$ by minimizing the loss function $\mathcal{L}_{sum}^{k}\left(\boldsymbol{\theta}_{u},\boldsymbol{\theta}_{\rho}\right)$ continuously $T_u$ times. Then fix $\boldsymbol{\theta}_u$ and update $\boldsymbol{\theta}_{\rho}$ by minimizing the loss function $\mathcal{L}_{sum}^{k}\left(\boldsymbol{\theta}_{u}, \boldsymbol{\theta}_{\rho}\right)$ continuously $T_\rho$ times. In an inner loop, $\mu_{k}$ and $\lambda_{k}$ do not change with respect to $k$. Let $R$ and $K$ be the numbers of the inner loop and the outer loop, respectively. The total number of iteration $Epoch=K \times R$. 

Algorithm $\ref{Alg1}$ demonstrates a general algorithmic framework of pressent adversarial neural network method for the eigenvalue optimization problem.

\begin{algorithm}[htbp]
\caption{Adversarial neural network method for topology optimization of eigenvalue problems}
\SetAlgoLined
\label{Alg1}
\DontPrintSemicolon 
\KwIn{Choose the proper loss function $\mathcal{L}_{sum}^{k}\left(\boldsymbol{\theta}_{u}, \boldsymbol{\theta}_{\rho}\right)$ according to the problem.}
\KwData{$K$: number of outer loop iterations; $R$: number of inner loop iterations; $T_u$: number of inner iterations for training $\hat{u}$; $T_\rho$:  number of inner iterations for training $\hat{\rho}$; $\eta_u$: learning rate of $\hat{u}$; $\eta_\rho$: learning rate of $\hat{\rho}$.}
\KwData{$\lambda^{0}=0$, $\mu^{0}$, $\beta$, $S$.}
\For{$k=1$ \KwTo $K$}{
  \For{$i=1$ \KwTo $R$}{
    \For{$t=1$ \KwTo $T_u$}{
      Compute the loss function $\mathcal{L}_{sum}^{k}\left(\boldsymbol{\theta}_{u}, \boldsymbol{\theta}_{\rho}\right)$.\;
      Compute gradient about $\boldsymbol{\theta}_{u}$ in  $\mathcal{L}_{sum}^{k}\left(\boldsymbol{\theta}_{u}, \boldsymbol{\theta}_{\rho}\right)$, and update $\boldsymbol{\theta}_{u}$ by Adam method.\;
    }
    \For{$t=1$ \KwTo $T_\rho$}{
      Compute the loss function $\mathcal{L}_{sum}^{k}\left(\boldsymbol{\theta}_{u}, \boldsymbol{\theta}_{\rho}\right)$.\;
      Compute gradient about $\boldsymbol{\theta}_{\rho}$ in $\mathcal{L}_{sum}^{k}\left(\boldsymbol{\theta}_{u}, \boldsymbol{\theta}_{\rho}\right)$, and update $\boldsymbol{\theta}_{\rho}$ by Adam method.\;
    }
  }
  $\mu_{k}=\min(\beta\mu_{k-1},S)$.\;
    $\lambda_{k}=\lambda_{k-1}+2 \mu_{k-1} h\left( \boldsymbol{\theta}_{\rho}\right)$ for augmented Lagrangian method.\;
}
\end{algorithm}

\section{Numerical results}
This section shows using adversarial neural network algorithm on \eqref{comp_question},
\eqref{comp_question_max}, \eqref{min1_four} and \eqref{min2_four}. 
All numerical examples were computed on the Google Colab platform.

\subsection{Second-order problems}
\textbf{Example 1: }Consider an annular domain $\Omega=\left\{(x, y): r_{\rm in}^{2}
\leq x^{2}+y^{2} \leq 1\right\}$. Let $\tau= 1 / r_{\text {in }}$ (indicating the ratio of the outer radius to the radius of the inner circle), $\delta=c /|\Omega|$ (indicating the proportion of the area constraint to the total area). Set a uniform grid of 50 $\times$ 50 in the square $[-1,1]^2$ and then choose the grid points in the annular as the sampling points. 

The boundary condition function is constructed by
$\ell_{u}=(1-x^{2}-y^{2})(x^{2}+y^{2}-r_{\rm in}^{2}).$
Set
\begin{equation}\label{u}
\hat{u}\left(x,y ; \boldsymbol{\theta}_{u}\right)=\ell_{u}(x,y) H\left(x,y ; \boldsymbol{\theta}_{u}\right),
\end{equation}
where $H\left(x,y ; \boldsymbol{\theta}_{u}\right)$ is the neural network.
For function $\rho\in [0,1]$, suitable estimation allows us to make the following construction of $\hat{\rho}$:
\begin{equation}\label{rho}
\hat\rho(x, y)=\max \left(0, \min \left(1, -g(x, y) H\left(x,y;\boldsymbol{\theta}_{\rho}\right)\right)\right)+1,
\end{equation}
with
$$
g(x, y)=\min(1-\sqrt{x^{2}+y^{2}},\sqrt{x^{2}+y^{2}}-r_{\rm in}).
$$

We apply the augmented Lagrangian method to treat the geometry constraint. 
Set Algorithm parameters: $T_u=4,T_\rho=1,R=50,K=15,{\eta}_{u}=10^{-4},{\eta}_{\rho}=10^{-3}$, $\beta=1.1,\mu_0=10,$ and $S=1000$. Set $\rho_1=0$ and $\rho_2=1$.

\begin{figure}[htbp]
\centering
\subfigure[$\alpha=1,\delta=0.64,\tau=3.5$]{
\includegraphics[width=5cm]{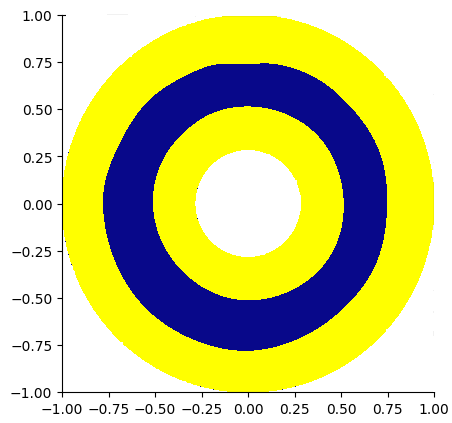}
}
\quad
\subfigure[$\alpha=1,\delta=0.64,\tau=1.2$]{
\includegraphics[width=5cm]{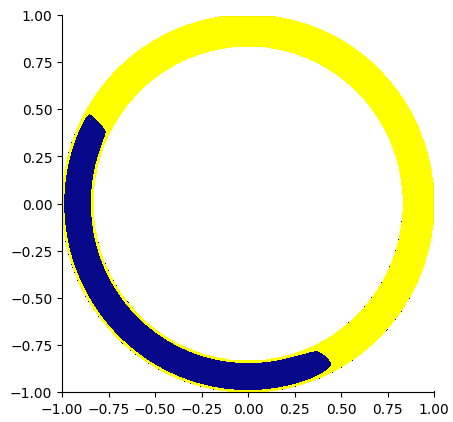}
}
\\
\subfigure[$\alpha=10,\delta=0.64,\tau=3.5$]{
\includegraphics[width=5cm]{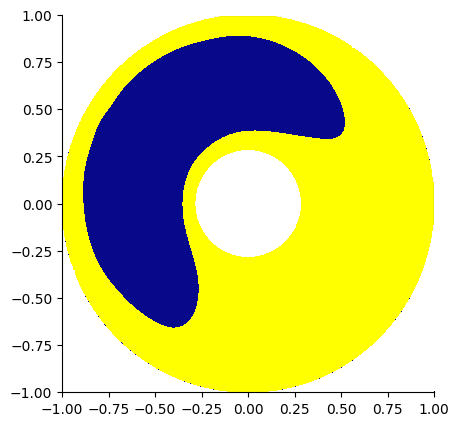}
}
\quad
\subfigure[$\alpha=1,\delta=0.83,\tau=3.5$]{
\includegraphics[width=5cm]{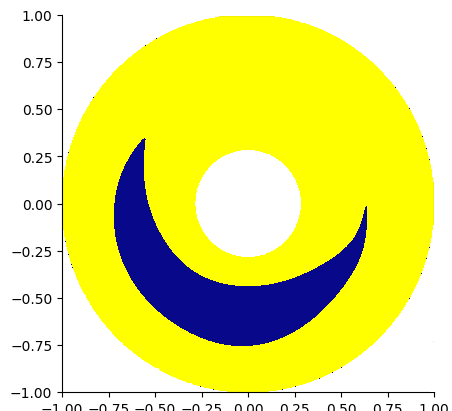}
}
\caption{Optimized designs for elliptic eigenvalue optimization for Example 1.}\label{yh_fig}
\end{figure}

\begin{figure}[htbp]
\centering
\subfigure[$\alpha=1,\delta=0.64,\tau=3.5$]{
\includegraphics[width=6cm]{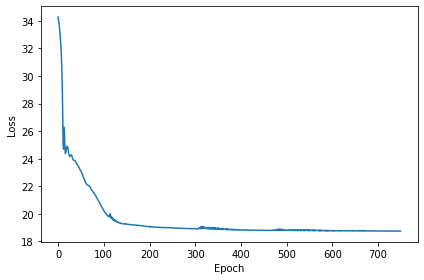}
}
\quad
\subfigure[$\alpha=1,\delta=0.64,\tau=1.2$]{
\includegraphics[width=6cm]{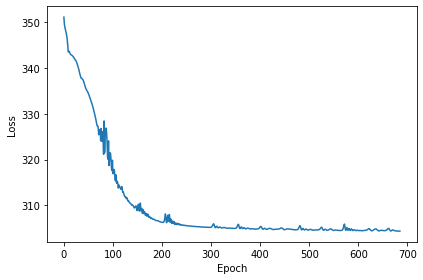}
}
\quad
\subfigure[$\alpha=10,\delta=0.64,\tau=3.5$]{
\includegraphics[width=6cm]{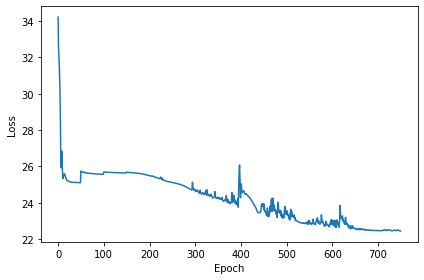}
}
\quad
\subfigure[$\alpha=1,\delta=0.83,\tau=3.5$]{
\includegraphics[width=6cm]{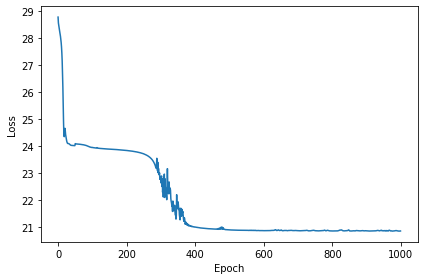}
}
\caption{Convergence of loss functions on elliptic eigenvalue optimization for Example 1.}\label{yh_epoch_fig}
\end{figure}

In Fig. \ref{yh_fig}, blue and yellow regions represent $\Omega\setminus D$ and yellow $D$, respectively. If $\alpha$ or $\delta$ is relatively not large, i.e., ($\alpha=1$ or $\delta=0.64$), an optimized design with radial symmetry is obtained as shown in Fig. \ref{yh_fig} (a). But in a same annular $\Omega$, Figs. \ref{yh_fig} (c)-(d) show radial symmetry breaking for final designs when increasing $\alpha$ from 1 to 10 or increasing the target area of $D$. Compared with Fig. \ref{yh_fig} (a), Fig. \ref{yh_fig} (b) also has radial symmetry breaks when reduce the width of the annulus. These observed phenomenons are consistent with that describes in \cite{CSG2000}. Fig. $\ref{yh_epoch_fig}$ shows that the objectives converge.

\begin{figure}[htbp] 
\centering 
\includegraphics[scale=0.4]{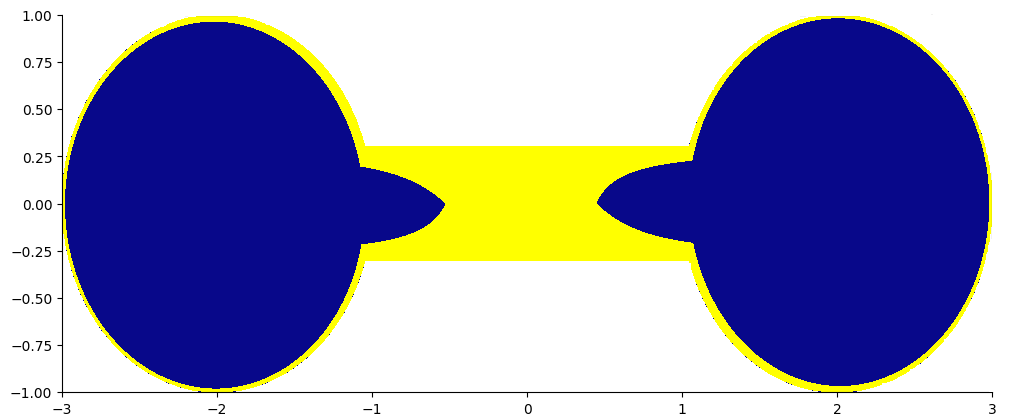} 
\caption{Optimized shape of $D$ in dumbbell domain for Example 2.}\label{yaling_fig}
\end{figure}

\textbf{Example 2: }Consider a dumbbell domain $
\Omega=B_{1}(-2,0) \cup((-2,2) \times(-0.3, 0.3)) \cup B_{1}(2,0),
$
where $B_{1}(p)=\left\{x \in \mathbb{R}^{2}:|x-p|<r\right\}$. We still use the uniform sampling method, establish a regular grid of 60 $\times$ 60 in the rectangle area [-3,3]$\times$[-1,1] which surrounds the dumbbell, and then take the square grid points in the dumbbell as the sampling point. There is a difference between the two adversarial neural networks $H\left(x,y ; \boldsymbol{\theta}_{u}\right)$, $H\left(x,y ; \boldsymbol{\theta}_{\rho}\right)$ we established in dumbbells domain with those in annular domain that in output layer, it is no longer simply linear but is square function. $\hat{u}$, $\hat{\rho}$ are bulit according to $\eqref{u}$, $\eqref{rho}$ respectively, in which $\ell_{u}$ is replaced as follows with symbolic distance function and dumbbell features.
$$
\begin{aligned}
\ell_{u}&=\max(0.3-\mid y \mid,\ell_{d})\ \chi_{\mid x \mid\leq 2.05}+\ell_{d}\ \chi_{\mid x \mid \geq 2.05},\\
\ell_{d}&=\max\left[1-(x-2)^{2}-y^{2},1-(x+2)^{2}-y^{2}\right].
\end{aligned}
$$
In the problem of dumbbell domain, $\alpha=0.1,\delta=0.3$ in \ref{comp_question}, we apply varying coefficient penalty method and obtain optimized shape of $D$ in Fig. \ref{yaling_fig}, which agrees well with that in \cite{CSG2000}. 

\begin{figure}[htbp]
\centering
\subfigure{
\includegraphics[width=4.6cm]{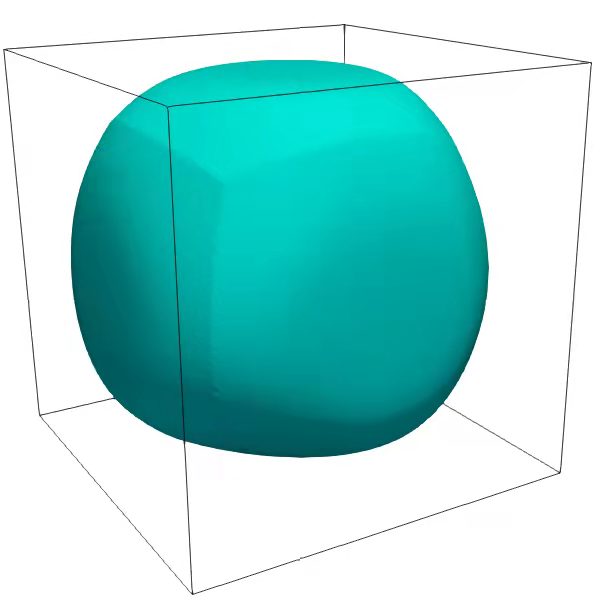}
}
\quad
\subfigure{
\includegraphics[width=7cm]{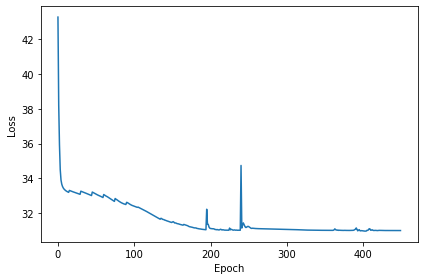}
}
\caption{Elliptic eigenvalue minimization in 3D: optimized design (left) and convergence history of loss (right) for Example 3.}\label{3d_fig}
\end{figure}

\textbf{Example 3: }For 3D, set $\alpha=10, c=0.5$. Let $\Omega=(0,1)^3$. We establish a regular grid of $40\times40\times40$ in the cube through uniform sampling. The adversarial neural networks here are $H\left(x,y,z ; \boldsymbol{\theta}_{u}\right)$ and $H\left(x,y,z ; \boldsymbol{\theta}_{\rho}\right)$.
We construct boundary condition function $\ell_{u}=64\left(xyz(1-x)(1-y)(1-z)\right)$ to satisfy boundary condition. Let
$\hat{u}\left(x,y,z ; \boldsymbol{\theta}_{u}\right)=\ell_{u}(x,y,z) H\left(x,y,z; \boldsymbol{\theta}_{u}\right),$
As to the minimization problem $\eqref{comp_question}$, we reasonably guess that $\rho=1$ on the boundary, $\rho=0$ at the central point. Choose
$$
\hat\rho(x, y,z)=\max \left(0, \min \left(1, \ell_{u}(x, y,z) \ell_{c}(x, y,z) H\left(x,y,z;\boldsymbol{\theta}_{\rho}\right)+g(x, y,z)\right)\right),
$$
where
$$
\begin{gathered}
\ell_{c}(x, y,z)=(x-0.5)^{2}+(y-0.5)^{2}+(z-0.5)^{2}, \quad
 g(x,y,z)=1-\ell_{u}^{'}(x, y,z)
\end{gathered}
$$
with  $\ell_{u}^{'}(x, y,z)=2\min(x,1-x,y,1-y,z,1-z)$.
Set $T_u=3$, $T_\rho=1$, $R=15$, $K=15$, ${\eta}_{u}= 10^{-4}$, and ${\eta}_{\rho}=5\times 10^{-4}$. We apply the augmented Lagrangian method for the geometry constraint and choose $\beta=1.5$, $\mu_0=1$, and $S=1000$. Finally, we obtain the optimized design $D$ and the convergence process of the loss function as shown in Fig. $\ref{3d_fig}$.

\begin{figure}[htbp]
\centering
\subfigure{
\includegraphics[width=4.6cm ]{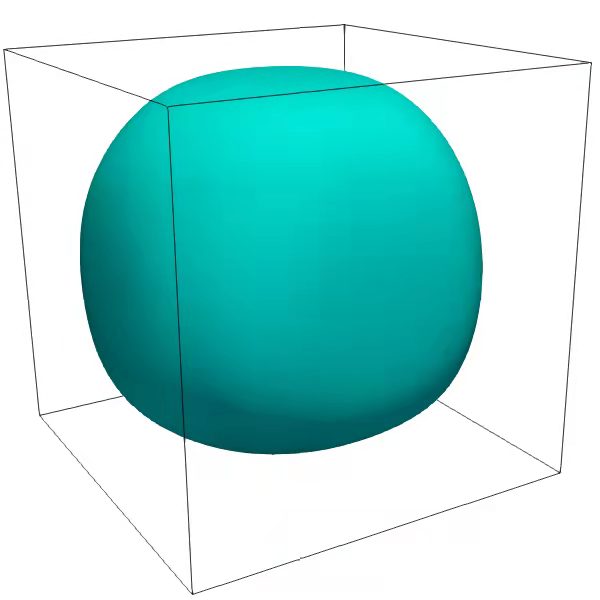}
}
\quad
\subfigure{
\includegraphics[width=7cm ]{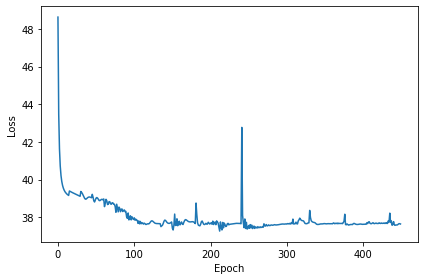}
}
\caption{Elliptic eigenvalue maximization in 3D: optimized design (left) and convergence history of loss (right) for Example 3.}\label{3d_fig_max}
\end{figure}

As to maximization problem $\eqref{comp_question_max}$, we reasonably set $\rho=0$ on the boundary and $\rho=1$ at the central point. We choose
$$
\hat\rho(x, y,z)=\max \left(0, \min \left(1, \ell_{u}(x, y,z) \ell_{c}(x, y,z) H\left(x,y,z;\boldsymbol{\theta}_{\rho}\right)+g(x, y,z)\right)\right),
$$
where
$$
\begin{gathered}
\ell_{c}(x, y,z)=(x-0.5)^{2}+(y-0.5)^{2}+(z-0.5)^{2}, \quad g(x,y,z)=\ell_{u}'(x, y,z).
\end{gathered}
$$
The final optimized design and the convergence process of the loss function are shown in Fig. $\ref{3d_fig_max}$.

In Fig. \ref{3d_fig}, the optimal design of minimization problem is the complementary set of solid sphere. The first eigenvalue $\lambda_1$ is optimized to a minimal value of 30.07.  To assess the accuracy of the maximization of the first eigenvalue, we employ the finite element method and utilize Freefem++ software \cite{ZZ2019}  to solve the corresponding eigenvalue problem within this region. The minimum eigenvalue obtained is 29.94, which is close to the result of neural network.
In Fig. \ref{3d_fig_max}, the optimal region $D$ of maximization problem in 3D is the solid sphere in the center. We get a maximal eigenvalue of 39.32 which is also close to the result of 39.20 computed by finite element method. Our model's results on optimization in 3D behave well as expected.

\subsection{Fourth-order problems}

\subsubsection{Eigenvalue optimization in inhomogeneous clamped plate}

For problem $\eqref{min1_four}$, choose $c=1.5$, $\rho_1=1$, and $\rho_2=2$. Let $\Omega$ be a square, circle and annulus. In particular, as to clamped plate problem, since it is easy to find boundary condition function satisfying boundary conditions, we use Full-satisfied boundary method. The neural networks that we build for $H\left(\boldsymbol{x}; \boldsymbol{\theta}_{u}\right)$ and $H\left(\boldsymbol{x}; \boldsymbol{\theta}_{\rho}\right)$ have 3 residual blocks. Every block has two fully connected layers and every layer has 80 neurons. All the activation functions are tanh except that we take linear function in the output layer. Set $T_u=3,T_\rho=1,R=15,K=30,{\eta}_{u}=10^{-4},{\eta}_{\rho}=10^{-3}$, $\beta=2,\mu_0=1$, and $S=1000$. In this fourth-order biharmonic eigenvalue optimization problem, neural network structures and parameters are the same while the constructions of $\hat{u}$ and $\hat{\rho}$ vary from the shape of $\Omega$ and the optimization target. We minimizing (resp. maximizing) the eigenvalue for Examples 4-6 (resp. Example 7). 

\begin{figure}[htbp]
\centering
\subfigure{
\includegraphics[width=4.8cm]{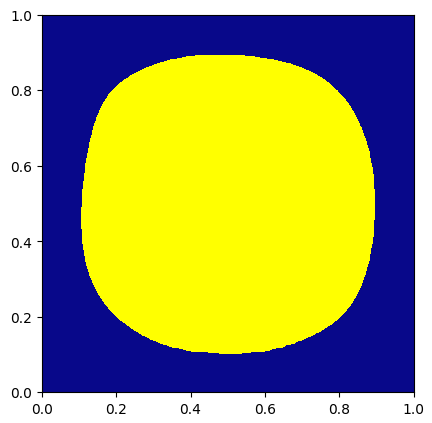}
}
\quad
\subfigure{
\includegraphics[width=6cm]{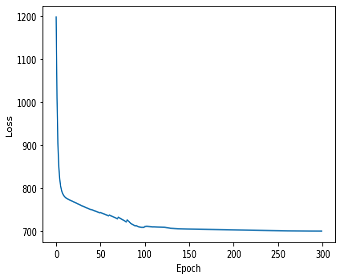}
}
\caption{Eigenvalue minimization of the inhomogeneous clamped plate for Example 4: optimized design (left) and convergence history of loss (right).}\label{clamped_sq_fig}
\end{figure}

\textbf{Example 4: }Consider $\Omega=(0,1)\times(0,1)$ and choose $\ell_{u}=256[x(x-1)y(y-1)]^{2}$ to
satisfy boundary conditions. As to $\hat{\rho}$, we reasonably guess that $\rho=1$ on the boundary and $\rho=2$ at the central point. Similar to the second-order problem, we construct
\begin{equation}\label{rho_4}
\hat\rho(x, y)=\max \left(1, \min \left(2, \ell_{u}(x, y) \ell_{c}(x, y) H\left(x,y;\boldsymbol{\theta}_{\rho}\right)+g(x, y)\right)\right),
\end{equation}
where
$$
\begin{gathered}
\ell_{c}(x, y)=(x-0.5)^{2}+(y-0.5)^{2}, \quad g(x, y)=\ell_{u}(x, y) +1.
\end{gathered}
$$
After training, the optimal region (figure $\ref{clamped_sq_fig}$) is a circle in the center while the blue region represents the density of 1 and the yellow region represents the density of 2. When $\Omega$ is a square domain, our method has been converged within 300 epochs and gets the minimal eigenvalue of 677.3 while the result obtained by the rearrangement algorithm \cite{CWC2016} is 652.2, having a 3.8\% difference. Each epoch takes 1.5 seconds.

\textbf{Example 5: }Consider $\Omega$ is a unit disk domain: $\Omega= \left\{(x, y): x^{2}+y^{2} \leq 1\right\}$, we choose $\ell_{u}=\left(1-x^2-y^2\right)^{2}$ to satisfy boundary conditions and then construct $\hat{u}$ by $\eqref{hat_u4}$. 
As to $\hat{\rho}$, we reasonably guess that $\rho=1$ on the boundary, $\rho=2$ at the central point and we construct 
\begin{equation}\label{rho_4_cy}
\hat\rho(x, y)=\max \left(1, \min \left(2, \ell_{u}(x, y) \ell_{c}(x, y) H\left(x,y;\boldsymbol{\theta}_{\rho}\right)+g(x, y)\right)\right),
\end{equation}
where
$$
\begin{gathered}
\ell_{c}(x, y)=x^{2}+y^{2}, \quad g(x, y)=\ell_{u}(x, y) +1.
\end{gathered}
$$

\begin{figure}[htbp]
\centering
\subfigure{
\includegraphics[width=5.3cm]{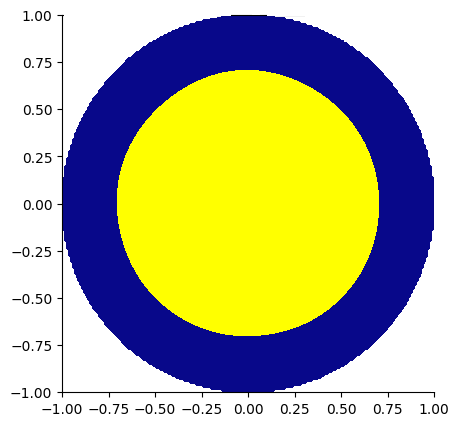}
}
\quad
\subfigure{
\includegraphics[width=6cm]{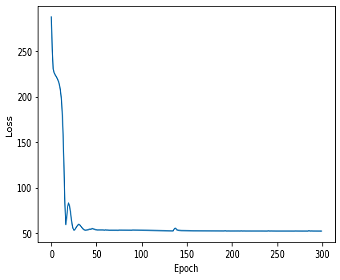}
}
\caption{Eigenvalue minimization of the inhomogeneous clamped plate for Example 5: optimized design (left) and convergence history of loss (right).}\label{clamped_cy_fig}
\end{figure}

After training, it is clear that our method has a good performance(figure $\ref{clamped_cy_fig}$), the optimal region $D$ is highly consistent with the shape optimized by the rearrangement algorithm \cite{CWC2016} and the minimal has converged to 52.3 within 50 epochs, which is close to 52.9 obtained by the rearrangement algorithm \cite{CWC2016}. 

\begin{figure}[htbp]
\centering
\subfigure{
\includegraphics[width=5.3cm]{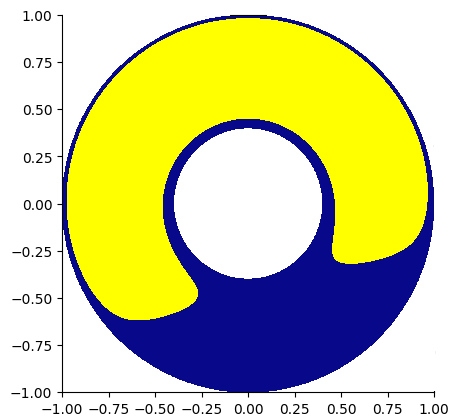}
}
\quad
\subfigure{
\includegraphics[width=6cm]{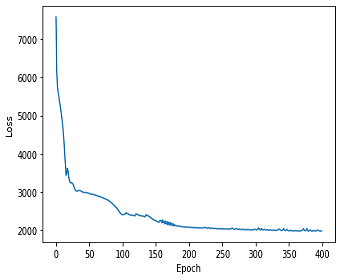}
}
\caption{Eigenvalue minimization of the inhomogeneous clamped plate for Example 6: optimized design (left) and convergence history of loss (right).}\label{clamped_yh_fig}
\end{figure}

\textbf{Example 6: }Consider an annulus $\Omega= \left\{(x, y): 0.4^2
\leq x^{2}+y^{2} \leq 1\right\}$ and choose $$\ell_{u}=\left(1-x^{2}-y^{2}\right)^2\left(x^{2}-y^{2}-0.16\right)^2.$$ We also guess that $\rho=1$ on the boundary and $\rho=2$ at the central to construct
$$
\hat\rho(x, y)=\max \left(1, \min \left(2, \ell_{u}(x, y) H\left(x,y;\boldsymbol{\theta}_{\rho}\right)\right)\right)+1,
$$
The first eigenvalue is optimized to the minimal of 1998 which differs the result of 1946 optimized by the rearrangement algorithm \cite{CWC2016} by 2.6\%. Fig. $\ref{clamped_yh_fig}$ shows convergence of loss and the optimized design agrees well with that in \cite{CWC2016}.

%

\begin{figure}[htbp]
\centering
\subfigure{
\includegraphics[width=5.3 cm]{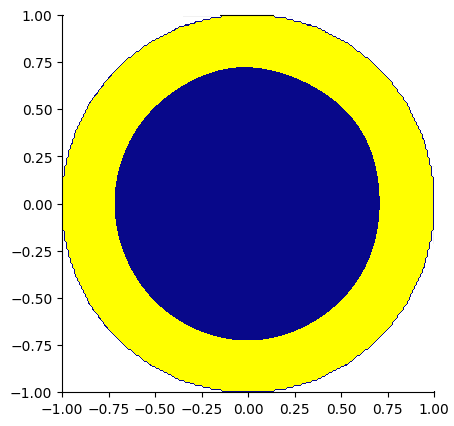}
}
\quad
\subfigure{
\includegraphics[width=6 cm]{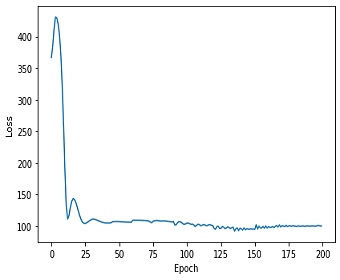}
}
\caption{Eigenvalue maximization of the inhomogeneous clamped plate for Example 7: optimized design (left) and convergence history of loss (right).}\label{clamped_cy_fig_max}
\end{figure}

\textbf{Example 7: }Consider $\Omega$ is a unit disk domain: $\Omega= \left\{(x, y): x^{2}+y^{2} \leq 1\right\}$, we construct $\hat{u}$ the same as that in the minimization problem. As to $\hat{\rho}$, contrary to the situation in minimization problem, we reasonably guess that $\rho=2$ on the boundary, $\rho=1$ at the central point and we construct it as follows,
\begin{equation}\label{rho_4_cy_max}
\hat\rho(x, y)=\max \left(1, \min \left(2, \ell_{u}(x, y) \ell_{c}(x, y) H\left(x,y;\boldsymbol{\theta}_{\rho}\right)+g(x, y)\right)\right),
\end{equation}
where 
$$
\begin{gathered}
\ell_{c}(x, y)=x^{2}+y^{2},\quad g(x, y)=\ell_{c}(x, y)+1.
\end{gathered}
$$
The maximized eigenvalue 101.8 obtained agrees well with that (a maximum of 101.7) in \cite{RSP2014} and the shape optimized is similar to that in \cite{RSP2014} (see Fig. $\ref{clamped_cy_fig_max}$ for design and history of loss).

\subsubsection{Eigenvalue optimization in inhomogeneous simply supported plate}
In the simply supported plate problem $\eqref{min2_four}$, we set $\nu=0.3$, $c=1.5$, $\rho_1=1$, and $\rho_2=2$. Considering that the boundary condition is relatively complex, we apply the half-satisfied boundary method \eqref{sum4_func_2} and \eqref{sum4_func_4}. The structure of  $H\left(\boldsymbol{x} ; \boldsymbol{\theta}_{u}\right)$ and $H\left(\boldsymbol{x} ; \boldsymbol{\theta}_{\rho}\right)$ are the same as that in clamped plate problem. We minimizing (resp. maximizing) the eigenvalue for Examples 8-9 (resp. Example 10).

The total loss function $\mathcal{L}_{sum}$ is $\eqref{sum4_func_2}$ and $\mathcal{L}_w$ in the $\eqref{sum4_func_2}$ is that of the augmented lagrangian method while $\mathcal{L}_{in}$ is given in $\eqref{sum4_func_2}$ 

\textbf{Example 8: }Consider a square $\Omega=(0,1)\times(0,1)$. We choose $\ell_{u}=16x(1-x)y(1-y)$. Then construct $\hat{u}$ and $\hat{\rho}$ according to $\eqref{hat_u4}$ and $\eqref{rho_4}$ respectively. In addition, $\kappa=0$ on the square boundary.
\begin{figure}[htbp]
\centering
\subfigure{
\includegraphics[width=4.8cm]{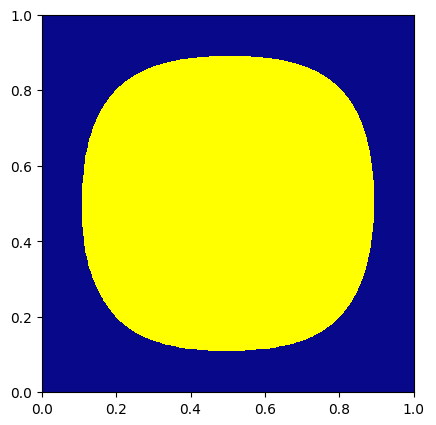}
}
\quad
\subfigure{
\includegraphics[width=4.9cm]{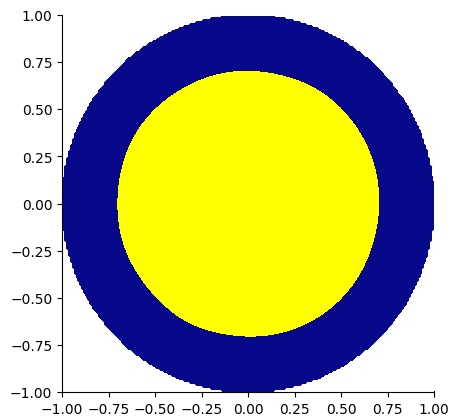}
}
\caption{Eigenvalue minimization of the inhomogeneous simply supported plate problem, optimal region $D$ when $\Omega$ is a square (left) for Example 8 and a unit disk (right) for Example 9.}\label{supported_sq_fig}
\end{figure}

\textbf{Example 9: }Consider a disk $\Omega= \left\{(x, y): x^{2}+y^{2} \leq 1\right\}$ and choose $\ell_{u}=1-x^2-y^2.$ Then construct $\hat{u}$ and $\hat{\rho}$ according to $\eqref{hat_u4}$ and $\eqref{rho_4_cy}$ respectively. Set $\kappa=1$ on the boundary In $\mathcal{L}_{in}$ and take $A_1=1$ and $A_2=2$.

\begin{figure}[htbp]
\centering
\subfigure{
\includegraphics[width=6.5cm]{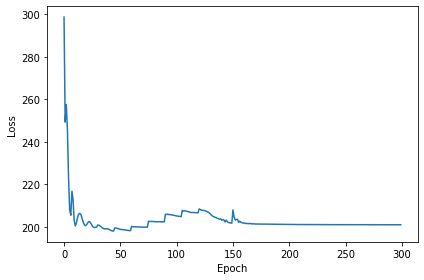}
}
\quad
\subfigure{
\includegraphics[width=6.5cm]{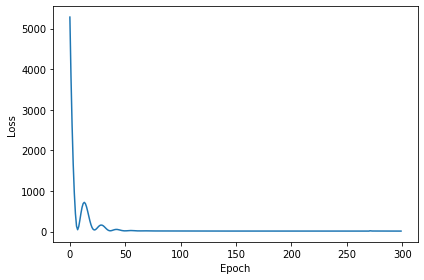}
}
\caption{Convergence history of loss for biharmonic eigenvalue minimization of the inhomogeneous simply supported plate: square for Example 8 (left) and disk for Example 9 (right).}\label{supported_cy_fig}
\end{figure}

It can be seen from figure $\ref{supported_sq_fig}$ and figure $\ref{supported_cy_fig}$ that the adversarial neural network method under the half-satisfied boundary method has good performance on both square and unit disk regions since the shape of the optimal density distribution under the two regions is consistent with that in the rearrangement algorithm \cite{CWC2016}. With our method, the first eigenvalue in the square domain is minimized to 202.04, the first eigenvalue in the unit disk domain is minimized to 12.83, while the result of the rearrangement algorithm is 201.60 and 12.84. so, our algorithm is still efficient and accurate.

\begin{figure}[htbp]
\centering
\subfigure{
\includegraphics[width=5.5cm]{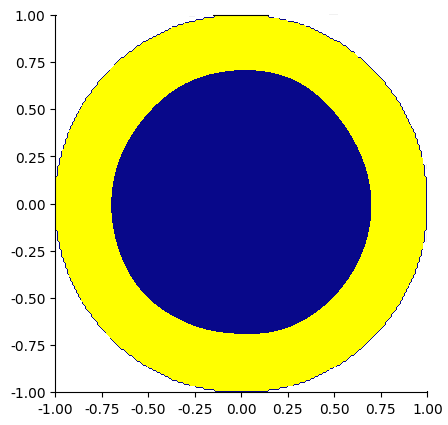}
}
\quad
\subfigure{
\includegraphics[width=6.8 cm]{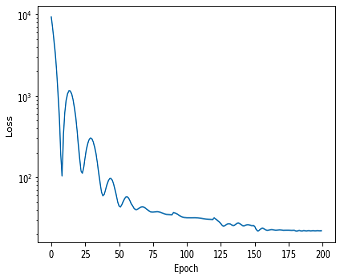}}
\caption{Eigenvalue maximization of the inhomogeneous simply supported plate for Example 10: optimized design (left) and convergence history of loss (right).}\label{supported_cy_fig_max}
\end{figure}

\textbf{Example 10: }Consider a unit disk $\Omega= \left\{(x, y): x^{2}+y^{2} \leq 1\right\}$ and choose $\ell_{u}=1-x^2-y^2$.
Then construct $\hat{u}$ and $\hat{\rho}$ according to $\eqref{hat_u4}$ and $\eqref{rho_4_cy_max}$ respectively. The total loss function $\mathcal{L}_{sum}$ is $\eqref{sum4_func_4}$ and loss $\mathcal{L}_w$ in $\eqref{sum4_func_4}$ is that of the augmented Lagrangian method while $\mathcal{L}_{in}$ is \eqref{obj42_func}. After training, the first eigenvalue of $ \lambda_1$ was finally optimized to a maximum of 22.3 and the optimal density distribution map is shown in Fig. \ref{supported_cy_fig_max}.

\section{Conclusions}

We have proposed an adversarial neural network topology optimization method based on deep learning to optimize second-order elliptic and fourth-order plate eigenvalues. The present algorithm does not need to solve eigenvalue problems repeatly, while it achieves the same optimized objective as the traditional sensitivity analysis-based methods which typically use partial differential equation solvers such as finite elements. Moreover, shape and topological changes can happen during training neural networks. Numerical examples are given for minimizing and maximizing the first eigenvalues in typical domains to show effectiveness and efficiency of the algorithm. In future, the adversarial neural network topology optimization method will be studied to solve more shape and topology optimization problems.


\begin{thebibliography}{99}
\bibitem{BAA2015}
B. Alipanahi, A. Delong, M.T. Weirauch, and B. J. Frey, Predicting the sequence specificities of DNA- and RNA-binding proteins by deep learning, Nat. Biotechnol. 33, 831–838, 2015.


\bibitem{BGY2020}
G. Bao, X. Ye, Y. Zang, and H. Zhou, Numerical solution of inverse problems by weak adversarial networks, Inverse. Probl. 36, 115003, 2020.

\bibitem{AGB2015}
A. G. Baydin, B. A. Pearlmutter, A. A. Radul, and J. M. Siskind, Automatic differentiation in machine learning: a survey, J. Machine. Learn. R. 18, 1-43, 2015.

\bibitem{MPB2003}
M. P. Bendsoe and O. Sigmund, Topology Optimization: Theory, Methods and Applications, Springer Press, Berlin, 2003.

\bibitem{JK2018}
J. Berg and K. Nyström. A unified deep artificial neural network approach to partial differential equations in complex geometries, Neurocomputing, 317, 28-41, 2018.

\bibitem{DPB2014}
D. P. Bertsekas, Constrained Optimization and Lagrange Multiplier Methods, Academic press, 2014.

\bibitem{TFC2003}
T. F. Chan and X. C. Tai, Level set and total variation regularization for elliptic inverse problems with discontinuous coefficients, J. Comput. Phys. 193, 40–66. 2003.

\bibitem{CSG2000}
S. Chanillo,  D. Grieser, M. Imai, K. Kurata, and I. Ohnishi, Symmetry breaking and other phenomena in the optimization of eigenvalues for composite membranes, Commun. Math. Phys. 214, 315–337, 2000.

\bibitem{CJSZ}
S. Cen, B. Jin, K. Shin, and Z. Zhou, Electrical impedance tomography with deep Calderón method, J. Comput. Phys. 493, 112427, 2023.

\bibitem{CWC2016}
W. Chen, C. S. Chou, and C. Y. Kao, Minimizing eigenvalues for inhomogeneous rods and plates, J. Sci. Comput. 69, 983-1013, 2016.

\bibitem{CSJ2000}
S. J. Cox and D. C. Dobson, Band structure optimization of two-dimensional photonic crystals in $H$-polarization. J. Comput. Phys. 158, 214–224, 2000.

\bibitem{Deng2021}
H. Deng and A. C. To,  A Parametric Level Set Method for Topology Optimization Based on Deep Neural Network.  ASME. J. Mech. Des. 143(9): 091702, 2021.




\bibitem{WEA2017}
W. E, A proposal on machine learning via dynamical systems, Comm. Meth. Stat. 5, 1-11, 2017.

\bibitem{WEB2018}
W. E and B. Yu, The Deep Ritz method: A deep learning based numerical algorithm for solving variational problems, Comm. Meth. Stat. 6, 1-12, 2018.

\bibitem{GXB2011}
X. Glorot, A. Bordes, and Y. Bengio,  Deep sparse rectifier neural networks, In Proceedings of the Fourteenth International Conference on Artificial Intelligence and Statistics, 315–323, 2011.

\bibitem{IGY2016}
I. Goodfellow, Y. Bengio, and A. Courville, Deep Learning, MIT Press, Cambridge, 2016.






\bibitem{KHX2016}
K. He, X. Zhang, S. Ren, and  J. Sun, Deep residual learning for image recognition, IEEE Conference on Computer Vision and Pattern Recognition (CVPR), 770-778, 2016.

\bibitem{ZZ2019}
F. Hecht, New development in Freefem++, J. Numer. Math. 20, 251-265, 2012.

\bibitem{Henrot}
A. Henrot, Extremum Problems for Eigenvalues of Elliptic Operators, Springer Science \& Business Media, 2006.

\bibitem{JCA2020}
J. Hendriks, C. Jidling, A. Wills, and T. Schön, Linearly constrained neural networks, arXiv preprint arXiv:2002.01600, 2020.

\bibitem{Hu2023}
Xindi Hu, Meizhi Qian and Shengfeng Zhu, Accelerating a phase field method by linearization for eigenfrequency topology optimization. Struct. Multidiscip. Optim. 66 (2023) 242. 

\bibitem{JHH2022}
J. Huang, H. Wang, and T. Zhou, An augmented Lagrangian deep learning method for variational problems with essential boundary conditions, Commu. Comput. Phys. 31, 966-986, 2022.

\bibitem{Huang2023}
M. Huang, Z. Du, C. Liu, Y. Zheng, T. Cui, Y. Mei, X. Li, X. Zhang, and X. Guo, A Problem-Independent Machine Learning (PIML) enhanced substructure-based approach for large-scale structural analysis and topology optimization of linear elastic structures, Extreme Mech. Lett. 63, 2023.

\bibitem{KK}
D. Kang and C.-Y. Kao, Minimization of inhomogeneous biharmonic eigenvalue problems. Appl. Math. Model. 51 (2017) 587-604.

\bibitem{DKJ2016}
D. Kingma and J. Ba, Adam: A method for stochastic optimization. Computer Science, 2014.

\bibitem{AKI2012}
A. Krizhevsky, I. Sutskever, and G. E. Hinton, Imagenet classification with deep convolutional neural networks, in: Advances in Neural Information Processing Systems, 1097–1105, 2012.

\bibitem{Lei2019}
X. Lei, C. Liu,  Z. Du, W. Zhang,  and X.  Guo, Machine Learning-Driven Real-Time Topology Optimization Under Moving Morphable Component-Based Framework, ASME. J. Appl. Mech.  86(1): 011004, 2019.


\bibitem{HPB2020}
B. M. Lake, R. Salakhutdinov, and J. B. Tenenbaum, Human-level concept learning through probabilistic program induction, Science 350, 1332–1338, 2015.

\bibitem{LLY}
K. Liang, X. Lu, and J. Z. Yang, Finite element approximation to the extremal eigenvalue problem for inhomogenous materials,
Numer. Math. 130, 741–762, 2015.

\bibitem{LYY2006}
Y. Lou and E. Yanagida, Minimization of the principal eigenvalue for an elliptic boundary value problem with indefinite weight, and applications to population dynamics, Jpn. J. Ind. Appl. Math. 23, 275–292, 2006.

\bibitem{LLP2020}
L. Lu, R. Pestourie, W. Yao, Z. Wang, F. Verdugo, and S. G. Johnson,  Physics-informed neural networks with hard constraints for inverse design, SIAM, J. Sci. Comput. 43, B1105-B1132, 2021.

\bibitem{LLX2020}
L. Lu, X. Meng, and Z. Mao. DeepXDE: A deep learning library for solving differential equations, SIAM Rev. 63, 208–228, 2021.

\bibitem{MB}
S.A. Mohammadi and F. Bahrami, Extremal principal eigenvalue of the bi-Laplacian operator, Appl. Math.
Model. 40, 2291-2300, 2016

\bibitem{Oh2019}
S. Oh,  Y. Jung,  S. Kim,  I. Lee,  and N. Kang,   Deep Generative Design: Integration of Topology Optimization and Generative Models.  ASME. J. Mech. Des.  141(11): 111405, 2019.

\bibitem{SOF2001}
S. Osher and F. Santosa, Level set methods for optimization problems involving geometry and constraints I. Frequencies of a two-density inhomogeneous drum, J. Comput. Phys. 171, 272–288, 2001.

\bibitem{SOR2002}
S. Osher and R. Fedkiw, Level Set Methods and Dynamic Implicit Surfaces, New York: Springer-Verlag, 2002.


\bibitem{RSP2014}
R. S. Pedro Antunes, Optimal bilaplacian eigenvalues, SIAM, J. Control Optim. 52, 2250–2260, 2014.

\bibitem{Qian2022}
M. Qian and X. Hu and S. Zhu, A phase field method based on multi-level correction for eigenvalue topology optimization. Comput. Methods Appl. Mech. Eng. 401 (2022) 115646. 



\bibitem{RMP2019}
M. Raissi, P. Perdikaris, and G. E. Karniadakis,  Physics-informed neural networks: A deep learning framework for solving forward and inverse problems involving nonlinear partial differential equations. J. Comput. Phys. 378, 686-707, 2019.

\bibitem{JSKS2018}
J. Sirignano and K. Spiliopoulos, DGM: A deep learning algorithm for solving partial differential equations, J. Comput. Phy. 375, 1339-1364, 2018.

\bibitem{JS1992}
J. Sokolowski and J. P. Zolesio, Introduction to the Shape Optimization: Shape Sensitivity Analysis, Springer, Heidelberg, 1992.


\bibitem{WX}
Y. Wang and H. Xie, Computing multi-eigenpairs of high-dimensional eigenvalue problems using tensor neural networks, J. Comput. Phys. 506 (2024), Paper No. 112928.

\bibitem{ZYB2020}
Y. Zang, G. Bao, X. Ye, and H. Zhou,  Weak adversarial networks for high-dimensional partial differential equations, J. Comput. Phys. 411, 109409, 2020.


\bibitem{ZS2010}
S. Zhu, C. Liu, and Q. Wu,  Binary level set methods for topology and shape optimization of a two-density inhomogeneous drum, Comput. Methods Appl. Mech. Eng. 199, 2970-2986, 2010.

\end{thebibliography}
\end{document}